\documentclass[fleqn]{mat01}
\usepackage{times,mathtimy,amssymb,latexsym,boxit}
\begin{document}

\setcounter{page}{65}
\firstpage{65}

\font\xxx=msam10 at 10pt
\def\ab{\mbox{\xxx{\char'245}}}

\font\www=msam10 at 10pt
\def\gtr{\mbox{\www{\char'77}}}

\font\xx=tib at 10pt
\def\N{\mbox{\xx{N}}}
\def\C{\mbox{\xx{C}}}
\def\P{\mbox{\xx{P}}}

\newcommand{\R}{\mathbb{R}}
\newcommand{\CA}{{\mathcal{A}}}

\newtheorem{theo}{Theorem}
\renewcommand\thetheo{\arabic{section}.\arabic{theo}}
\newtheorem{theor}[theo]{\bf Theorem}
\newtheorem{lem}[theo]{Lemma}
\newtheorem{pot}[theo]{Proof of Theorem}
\newtheorem{propo}{\rm PROPOSITION}
\newtheorem{rema}[theo]{Remark}
\newtheorem{defn}[theo]{\rm DEFINITION}
\newtheorem{exam}{Example}
\newtheorem{coro}[theo]{\rm COROLLARY}
\def\conjecture{\trivlist\item[\hskip\labelsep{\it Conjecture.}]}

\newtheorem{theore}{\bf Theorem}
\renewcommand\thetheore{\arabic{theore}}
\newtheorem{case}[theore]{Case}

\def\suffici{\trivlist\item[\hskip\labelsep{{\it Sufficiency.}}]}

\renewcommand{\theequation}{\thesection\arabic{equation}}

\title{On the limit-classifications of even and odd-order formally
symmetric differential expressions}

\markboth{K~V~Alice, V~Krishna Kumar and A~Padmanabhan}{On the
limit-classifications of differential expressions}

\author{K~V~ALICE, V~KRISHNA KUMAR$^{*}$ and A~PADMANABHAN$^{\dagger}$}

\address{Department of Mathematics, Newman College, Thodupuzha~685~585,
India\\
\noindent $^{*}$Department of Mathematics, University of Calicut 673~635, India\\
\noindent $^{\dagger}$Department of Mathematics, Govt. College,
Mokeri~673~516, India\\
\noindent E-mail: vellatkrishna@yahoo.com}

\volume{114}

\mon{February}

\parts{1}

\Date{MS received 10 July 2003}

\begin{abstract}
In this paper we consider the formally symmetric differential expression
$M[\cdot]$ of any order (odd or even) $\geq 2$. We characterise the
dimension of the quotient space $D(T_{\max})/D(T_{\min})$ associated
with $M[\cdot]$ in terms of the behaviour of the determinants
\begin{equation*}
\det\limits_{r,s\in {\bf N}_{n}} [[f_{r}g_{s}](\infty)]
\end{equation*}
where $1\leq n\leq$ (order of the expression $+1$); here $[fg](\infty) =
\lim\limits_{x\rightarrow\infty}[fg](x)$, where $[fg](x)$ is the
sesquilinear form in $f$ and $g$ associated with $M$. These results
generalise the well-known theorem that $M$ is in the limit-point case at
$\infty$ if and only if $[fg](\infty) = 0$ for every $f,g\in$ the
maximal domain $\Delta$ associated with $M$.
\end{abstract}


\maketitle

\section{Introduction}

Let $\N$ denote the set of natural numbers and ${\N}_{k}: =
\{1,2,\ldots,k\}$ for $k\in\N$. We write $C^{(r)}(I)\ (r =
0,1,2,\ldots,m)$ for the class of complex valued functions defined on
the interval $I$ with $r$ continuous derivatives, and $AC_{\rm loc}(I)$ for
the functions which are absolutely continuous on all compact sub
intervals of $I$.

We consider the formally symmetric differential expression $M$ of order
$m\ (m=2k\ \hbox{or}\ 2k-1, k=1,2,\ldots)$ given by (i.e. $M = M^{+}$, the
formal Lagrange adjoint of $M$)
\begin{equation}
M[y] = \sum\limits_{r=0}^{h} (-1)^{r} (s_{r}y^{(r)})^{(r)} + \frac{1}{2}
\sum\limits_{r=0}^{k-1} i^{2r+1} \{(q_{r}y^{(r)})^{(r+1)} +
(q_{r}y^{(r+1)})^{(r)}\},
\end{equation}
where the sets of coefficients $\{s_{r}\}$ and $\{q_{r}\}$ are real
valued on $I$ with $I$ designating the semi infinite interval
$[0,\infty)$. Further we assume that
\begin{align}
s_{r} &\in C^{(r)}(I)\qquad (r = 0,1,2,\ldots, h)\nonumber\\
q_{r} &\in C^{(r+1)}(I)\qquad (r = 0,1,2,\ldots,k).
\end{align}
If $M$ is of even-order, $m = 2k\ (k=1,2,\ldots)$, then
\begin{equation}
h =k \ \hbox{in (1.1)}\quad \hbox{and}\quad s_{k}(x) > 0,\ x\in I.
\end{equation}
If $M$ is of odd-order, $m = 2k - 1\ (k = 1,2,\ldots)$, then
\begin{equation}
h = k - 1\ \hbox{in (1.1)}\quad \hbox{and}\quad q_{k-1}(x) > 0,\ x\in I.
\end{equation}
Indeed, any formally symmetric differential expression of order $m$ with
sufficiently smooth coefficients can be expressed as in (1.1) with
suitable choice of coefficient functions.

The conditions (1.2) and (1.3) or (1.4) show that $M$ is regular on
$[0,\infty)$, but $M$ has a singular point at $\infty$ (see \cite{12},
\S15.1).

The differential equation we are concerned with is given by
\begin{equation}
M[y] = \lambda y\quad \hbox{on}\ I,
\end{equation}
where $\lambda$ is a complex parameter $\lambda = \mu + i\nu$.

The standard existence theorems for ordinary, linear, homogeneous
differential equations apply to eq.~(1.5) (see the books: \cite{2}, ch.~3,
\S6 and \cite{12}, \S16.2).

The underlying Hilbert space for the analysis of the problem is the
collection of all equivalence classes of complex valued Lebesgue
measurable functions $f$ on $[0,\infty)$ such that
$\int_{0}^{\infty}|f|^{2} < \infty$ which is denoted by
$L^{2}(0,\infty)$.

The Green's formula for $M$ takes the form
\begin{equation}
\int_{0}^{x} \{\overline{g}M[f] - f\overline{M}[g]\} = [fg](x) -
[fg](0),\quad x\in (0,\infty)
\end{equation}
where $f^{(m-1)},g^{(m-1)}\in AC_{\rm loc}(0,\infty)$. Here the
integrated term $[fg](x)$ on the R.H.S. is a skew-hermition,
non-singular form on $[0,\infty)$ (see \cite{2}, ch.~3, \S6 and
\cite{8}, \S5).

Also
\begin{equation*}
[gf](x) = - [\overline{fg}](x),\quad x\in (0,\infty).
\end{equation*}

To set up the differential operators associated with $M[\cdot]$ in
$L^{2}(0,\infty)$, we introduce the linear manifold $\Delta$ defined by
\begin{equation*}
\Delta : = \{f: f\in L^{2} (0,\infty),\ \ f^{(m-1)} \in AC_{\rm
loc}(0,\infty)\ \ \hbox{and}\ \ M[f] \in L^{2}(0,\infty)\}.
\end{equation*}

$\left.\right.$\vspace{-1.5pc}

\noindent From an application of Green's formula we have
\begin{equation*}
[fg](\infty) := \lim\limits_{x\rightarrow \infty} [fg](x)
\end{equation*}
exists and is finite for all $f$ and $g$ in $\Delta$.

Next we introduce two differential operators $T_{\max}$ and $T_{\min}$,
associated with $M$, defined as follows:

\begin{enumerate}
\renewcommand{\labelenumi}{\arabic{enumi}.}
\item The maximal operator $T_{\max}$: The domain $D(T_{\max})$ is
$\Delta$ and $T_{\max}f = M[f]\ (f\in D(T_{\max}))$;
\item The minimal operator $T_{\min}$: The domain $D(T_{\min})$ is
\begin{equation*}
\hskip -1.25pc D(T_{\min}) = \{f: f\in \Delta\ \hbox{and}\ [fg](0) = [fg](\infty) = 0\
\hbox{for all}\ g \in \Delta\}
\end{equation*}
and
\begin{equation*}
\hskip -1.25pc T_{\min}f = M[f]\quad (f\in D(T_{\min})).
\end{equation*}
\end{enumerate}

The domains $D(T_{\min})$ and $D(T_{\max})$ are dense in
$L^{2}(0,\infty)$. These operators have the following properties:

\begin{enumerate}
\renewcommand{\labelenumi}{(\alph{enumi})}
\item $T_{\min}$ is a closed, symmetric operator in $L^{2}(0,\infty)$.
\item $T_{\max}$ is a closed, but not symmetric operator in
$L^{2}(0,\infty)$.
\item $T^{*}_{\min} = T_{\max}$; where $T_{\min}^{*}$ denotes the
adjoint operator of $T_{\min}$ (see \cite{3}, ch.~XIII, \S2.1--2.8 and
\cite{12}).
\end{enumerate}
The deficiency indices $(N_{+},N_{-})$ of the closed, symmetric operator
$T_{\min}$ are defined as
\begin{equation*}
N_{\pm} = \dim \{f:f\in D(T^{*}_{\min});\ T_{\min}^{*} f = \pm if\}.
\end{equation*}
From the general theory of deficiency indices of symmetric operators
(see \cite{3}, ch.~12, \S19), we have
\begin{align*}
N_{\pm} &= \dim \{y \in C^{m}[0,\infty); M[y] = \lambda y \ \hbox{on}\
[0,\infty),\\
&\qquad\quad\ y \in L^{2} (0,\infty)\ \hbox{and}\ \lambda \in \C_{\pm}\},
\end{align*}
where $C_{\pm} = \{\lambda \in \C,\ \hbox{Im}\ \lambda\ \gtr\ \ \ 0\}$.

Thus the deficiency indices $N_{+}(N_{-})$ represent the number of
linearly independent solutions of the differential equation~(1.5) which
are in $L^{2}(0,\infty)$ when $\lambda \in \C_{+} (\C_{-})$. Hence both
$N_{+}$ and $N_{-}$ are finite and
\begin{equation}
0\leq N_{\pm} \leq m,
\end{equation}
where $m$ is the order of the equation.

Further from the general theory of symmetric operators, it is known that
$T_{\min}$ has self-adjoint extensions, i.e. $T_{\max}$ has self-adjoint
restrictions in $L^{2}(0,\infty)$, if and only if $N_{+} = N_{-}$. A
better estimate of the lower bound for these indices in (1.7) are

\begin{enumerate}
\renewcommand{\labelenumi}{\arabic{enumi}.}
\item When $m = 2k\ (k=1,2,\ldots)$ we get
\begin{equation}
\hskip -1.25pc k\leq N_{+} \leq 2k = m,\ \ k\leq N_{-} \leq 2k = m.
\end{equation}

\item When $m = 2k - 1\ (k = 1,2,\ldots)$ we get
\begin{equation}
\hskip -1.25pc k-1\leq N_{+} \leq 2k-1 = m,\ \ k\leq N_{-} \leq 2k - 1 = m
\end{equation}
(see \cite{7} and \cite{8}).
\end{enumerate}

Any choice of integers $N_{+}(N_{-})$ satisfying (1.8) or (1.9) are
possible pairs of deficiency indices provided they are also subject to
the additional constraints $N_{+} = m$ if and only if $N_{-} = m$ for
$m\geq 2$. For $m = 1,\ N_{+} = 0$ and $N_{-} = 1$.

The differential expression $M$ is said to be in the limit $(N_{+},N_{-
})$ case at the singular point $\infty$ if the deficiency indices of the
corresponding minimal closed operator $T_{\min}$ in $L^{2}(0,\infty)$
are $(N_{+},N_{-})$. In particular, borrowing the terminology of Weyl
(see \cite{2}) we say that $M$ is in the limit-point case at $\infty$, if
$N_{+} = N_{-} = k$ in the even-order case and $N_{+} = k-1, N_{-} = k$
in the odd-order case.

We now introduce the Titchmarsh--Weyl $L^{2}(0,\infty)$ solutions of
(1.5) for the even- and odd-order cases.

\subsection*{\it The even-order case}

We assume $m = 2k\ (k\in \N)$. Let $\theta_{r}$ and $\phi_{r}\ (r\in
\N_{k})$ be solutions of (1.5) taking initial values at 0 which are
independent of $\lambda$, such that
\begin{equation*}
[\theta_{r}\theta_{s}] (0) = 0 = [\phi_{r}\phi_{s}](0),\quad
[\theta_{r}\phi_{s}](0) = \delta_{rs},\quad r,s\in \N_{k},
\end{equation*}
where $\delta_{rs}$ is the Kronecker delta function. Such a choice of
initial conditions is possible and the set $\{\theta_{r},\phi_{r};\
r\in\N_{k}\}$ forms a basis of solutions for (1.5). Then it can be seen
that there are $k^{2}$ analytic functions $\{m_{rs}(\cdot);\ r,s\in
\N_{k}\}$ which are all regular on $\C_{+}\cup \C_{-}$ and such that the
$k$ linearly independent solutions determined by
\begin{equation*}
\psi_{r}(x;\lambda): = \theta_{r}(x;\lambda) + \sum\limits_{s=1}^{k}
m_{rs} (\lambda)\phi_{s}(x;\lambda),\quad (x\in[0,\infty),\lambda \in
\C_{+}\cup \C_{-})
\end{equation*}

$\left.\right.$\vspace{-1.5pc}

\noindent belong to $L^{2}(0,\infty)$ for $r\in\N_{k}$ and for all $\lambda \in
\C_{+} \cup \C_{-}$. The analytic functions $m_{rs}(\lambda)$ satisfy
\begin{equation*}
m_{rs}(\lambda) = \overline{m}_{sr} (\overline{\lambda}),\quad
(r,s\in\N_{k}).
\end{equation*}

Thus it follows that to each $\lambda \in \C_{+} \cup \C_{-}$, there
exist at least $k$ linearly independent solutions of (1.5) which belong
to $L^{2}(0,\infty)$. The deficiency indices of the associated
$T_{\min}$ can be characterised as
\begin{equation*}
N_{\pm} = k +\,\hbox{\!dim of the}\ L^{2}(0,\infty)\ \hbox{span of}\
\{\phi_{r}(\cdot,\pm i), r\in \N_{k}\}.
\end{equation*}
For details, see \cite{7}.

\subsection*{\it The odd-order case}

Here $m = 2k - 1, k = 2,3,\ldots$. Let $\theta_{r}(x;\lambda)\ (x\in
[0,\infty),\lambda \in \C, r \in \N_{k-1})$ and $\phi_{s}(x,\lambda)\
(x\in [0,\infty), \lambda \in \C, s\in \N_{k})$ be solutions of (1.5)
taking initial values, independent of $\lambda$, at 0 such that
\begin{align}
[\theta_{r}\theta_{s}] (0) &= 0,\quad (r,s\in \N_{k-1})\nonumber\\[.2pc]
[\theta_{r}\phi_{s}] (0) &= \delta_{rs},\quad (r\in \N_{k-1},s\in
\N_{k})\nonumber\\[.2pc]
[\phi_{r}\phi_{s}] (0) &= i\delta_{rk}\delta_{sk},\quad (r,s\in\N_{k}).
\end{align}
Then the set of functions $\{\theta_{r},\phi_{s}; r\in \N_{k-1}, s\in
\N_{k}\}$ forms a basis for all solutions of (1.5). Further there exists
$k(2k - 1)$ analytic functions $\{p_{rs}(\cdot);r,s\in\N_{k}\}$ and
$\{n_{rs}(\cdot);\ r\in\N_{k},s\in\N_{k-1}\}$ with $p_{rs}(n_{rs})$
regular in $\C_{+}(\C_{-})$ such that if the solutions $\psi_{p,r}$ and
$\psi_{n,r}$ are defined by
\begin{align*}
\psi_{p,r}(x,\lambda) &= \theta_{r} (x,\lambda) + \sum\limits_{s= 1}^{k}
p_{rs}(\lambda) \phi_{s}(x,\lambda),\ (x\!\in\! [0,\infty), \lambda\! \in\!
\C_{+}, r\!\in\! \N_{k-1})\\[.2pc]
\psi_{p,k}(x,\lambda) &= \sum\limits_{s=1}^{k} p_{ks} (\lambda) \phi_{s}
(x,\lambda),\ \ (x\in [0,\infty), \lambda \in \C_{+})
\end{align*}
and
\begin{align}
\psi_{n,r}(x,\lambda) &= \theta_{r} (x,\lambda) + \sum\limits_{s= 1}^{k-1}
n_{rs}(\lambda) \phi_{s}(x,\lambda),\ (x\!\in\! [0,\infty), \lambda \!\in\!
\C_{-}, r\!\in\! \N_{k-1})\nonumber\\[.2pc]
\psi_{n,k}(x,\lambda) &= \phi_{k} (x,\lambda) + \sum\limits_{s=1}^{k-1}
n_{ks} (\lambda)\phi_{s}(x,\lambda),\ \ (x\in [0,\infty),\lambda \in
\C_{-})
\end{align}
then  $\psi_{p,r}(\cdot,\lambda) \in L^{2} (0,\infty)$ for $r\in
\N_{k},\ \lambda\in\C_{+}$ and $\psi_{n,r}(\cdot,\lambda)\in
L^{2}(0,\infty)$ for $r\in\N_{k},\lambda \in\C_{-}$ with the possibility
of $\psi_{p.k}(\cdot,\lambda)$ being a null solution of (1.5) in certain
cases (see \cite{10}, \S\S2 and 3). The connection between the existence of
the integrable square solutions and the deficiency indices of the
associated $T_{\min}$ is that
\begin{align}
N_{+} &= k-1+\ \dim L^{2}(0,\infty)\ \hbox{span of}\
\{\phi_{r}(\cdot,\lambda),r\in \N_{k},\lambda \in \C_{+}\},\nonumber\\
N_{-} &= k + \ \dim L^{2}(0,\infty)\ \hbox{span of}\
\{\phi_{r}(\cdot,\lambda), r\in \N_{k-1}, \lambda\in \C_{-}\}.
\end{align}
In both the even- and odd-order cases we have an elegant
characterisation of the limit-point cases in terms of the behaviour of
the sesquilinear form $[fg](x); f, g \in\Delta$ as $x\rightarrow$ the
singular point. We shall recall this result in the following theorem.

\begin{theor}[\!]
Let $M$ be a formally symmetric differential expression of order $m\ (m =
2k\ \hbox{or}\ 2k-1, k = 1,2,\ldots)$ given by $(1.1)$ on $[0,\infty)$ and
the coefficients satisfy the conditions $(1.2)$. A necessary and
sufficient condition for $M$ to be in the limit-point case {\rm (}i.e. limit
$(k,k)$ in the $2k${\rm th}-order case{\rm ,} limit $(k-1,k)$ in the $(2k-
1)${\rm th}-order case{\rm )} is that
\begin{equation*}
[fg](\infty) = \lim\limits_{x\rightarrow \infty} [fg](x) = 0\quad
\hbox{for all}\ f\ \hbox{and}\ g\ \hbox{in}\ \Delta.
\end{equation*}
{\rm (}See {\rm \cite{6}} and {\rm \cite{10}.)}
\end{theor}

The object of this paper is to generalise this
result. The generalisations are given by Theorems~2.1 and 2.2 for the
even- and odd-order cases respectively. Before we state these
generalisations in Theorems~2.1 and 2.2, we quote some known results
which find repeated application in the proof of the main results.

\setcounter{theo}{0}
\begin{lem}\hskip -.4pc {\rm (}A determinantal identity associated with
$M${\rm )}.\ \ Let $\{f_{r},g_{r}$ for $r = 1,2,\ldots,m+1\}$ be any two
sets of $(m+1)$ functions all in $\C^{(m)}[0,\infty)$. Then
\begin{equation*}
\det\limits_{r,s\in {\bf N}_{m+1}} [[f_{r}g_{s}] (x)] = 0,\quad
(x\in[0,\infty)),
\end{equation*}
where $[fg](\cdot)$ is the sesquilinear form in $f$ and $g$ associated
with $M$ {\rm (}see {\rm \cite{4}, \S11)}.
\end{lem}

For the description of the system $\{f_{r},g_{r}\}$ we use the following
convention:
\begin{align*}
f_{r}: f_{1},f_{2},\ldots, f_{m+1},\\
g_{r}: g_{1},g_{2},\ldots, g_{m+1},
\end{align*}
where $f_{r},g_{r}$ correspond to the specific functions to be
substituted in the identity.

\begin{lem}\hskip -.4pc {\rm (}the $L^{2}(0,\infty)$ lemma{\rm )}.\ \
Suppose that the complex valued measurable functions $f$ and $g$ on
$[0,\infty)$ are such that $f\in L^{2}(0,\infty),g\in L^{2}(0,x),(x\in
[0,\infty))$ and $g\notin L^{2}(0,\infty)$. Then
\begin{equation*}
\lim\limits_{x\rightarrow\infty} \left\lbrace \int_{0}^{x}
f\overline{g}\right\rbrace \left\lbrace
\int_{0}^{x}|g|^{2}\right\rbrace^{-1/2} = 0.
\end{equation*}
{\rm (}See {\rm \cite{6}, \S2.)}
\end{lem}

\section{The main results}

The main results of this paper are as follows:

\setcounter{theo}{0}
\begin{theor}[\!]
Even-order case. Let $M$ be a formally symmetric differential expression
of order $2k\ (k = 1,2,\ldots)$. Let $p$ and $q$ be non-negative integers
such that $0 \leq p,q\leq k$. Let $n = p+q$. Then a necessary condition
for $M$ to be in the limit $(k + p, k + q)$ case at $\infty$ is that
\begin{equation}
\det\limits_{r,s\in {\bf N}_{n+1}} [[f_{r}g_{s}](\infty)] =
0\quad \hbox{for all}\ f_{r},g_{s} \in \Delta.
\end{equation}
Conversely{\rm ,} if
\begin{equation}
\det\limits_{r,s\in {\bf N}_{n+1}} [[f_{r}g_{s}](\infty)] =
0\quad \hbox{for all}\ f_{r},g_{s}\in \Delta,
\end{equation}
then $M$ is in the limit $(k+p,k+q)$ case at $\infty$ where $p+q\leq n$.
More precisely $p+q = n$ if {\rm (2.14)} holds for all $f_{r},g_{s}\in \Delta
(r,s\in N_{n+1})$ but
\begin{equation*}
\det\limits_{r,s\in {\bf N}_{n}} [[f_{r}g_{s}](\infty)] \neq
0\quad \hbox{for some}\ f_{r},g_{s}\in \Delta.
\end{equation*}
\end{theor}

\setcounter{theo}{0}
\begin{rema}{\rm
Here the choice of $p$ and $q$ are not unique. Any $p,q$ such that $0
\leq p,q\leq k$ is possible subject to the other constraints. For
instance, in the case where $M$ is real, the deficiency indices are
necessarily equal and hence $p = q$. Also whether $M$ is real or not,
whenever $p = k$, then $q = k$. Such constraints do not come out from
the theorem.}
\end{rema}

\begin{theor}[\!]
Odd-order case. Let $M$ be a formally symmetric differential expression
of order $2k - 1\ (k = 2,3,\ldots)$. Let $p,q$ be non-negative integers
such that $0 \leq p \leq k,\ 0\leq q\leq k - 1$. Let $n = p + q$. Then a
necessary condition for $M$ to be in the limit $(k-1+p,k+q)$ case at
$\infty$ is that
\begin{equation}
\det\limits_{r,s\in {\bf N}_{n+1}} [[f_{r}g_{s}](\infty)] = 0\quad
\hbox{for all}\ f_{r},g_{s}\in \Delta.
\end{equation}
Conversely{\rm ,} if
\begin{equation}
\det\limits_{r,s\in {\bf N}_{n+1}} [[f_{r}g_{s}](\infty)] = 0\quad
\hbox{for all}\ f_{r},g_{s} \in \Delta,
\end{equation}
then $M$ is in the limit $(k-1 + p,k + q)$ case at $\infty$ for some $p$
and $q$ such that $p + q\leq n$. More precisely $p+q = n$ if {\rm (2.16)}
holds for all $f_{r},g_{s}\in\Delta\ (r,s\in\N_{n+1})$ but
\begin{equation*}
\det\limits_{r,s\in {\bf N}_{n}} [[f_{r}g_{s}](\infty)] \neq 0\quad
\hbox{for some}\ f_{r},g_{s}\in \Delta.
\end{equation*}
\end{theor}

\section{Proof of the results}

\setcounter{section}{2}
\setcounter{theo}{0}
\begin{pot}{\rm This is similar to the proof of Theorem~2.2, but less
complicated. Therefore we omit the details (see \cite{1}, ch.~4).

We choose to give the proof of the theorem for the odd-order case in
detail (see \cite{13}, ch.~5).}
\end{pot}

\begin{pot}\hskip -.4pc {\rm (}{\it Necessity}{\rm )}.\ \ {\rm Suppose that $M$ is in
the limit $(k-1 + p, k+q)$ case at $\infty$; where $p+q = n$. Let
$\lambda$ be a fixed point in $\C_{+}$ and $\phi(\cdot)$ and
$\tilde{\phi}(\cdot)$ denote $\phi(\cdot,\lambda)$ and
$\phi(\cdot,\overline{\lambda})$ respectively. Then, since $M$ is
assumed to be in the limit $(k-1+p,k+q)$ case at $\infty$, from (1.12)
it is clear that there exist numbers $t_{1},t_{2},\ldots, t_{p} \in
\N_{k}, t_{i} \neq t_{j}\ (i,j\in \N_{p})$ such that the
$L^{2}(0,\infty)$ spans of the $(k-p)$ functions $\{\phi_{s}(\cdot),
s\neq t_{i}, i\in \N_{p}\}$ are null. Also there exist numbers
$r_{1},r_{2},\ldots, r_{q}$ belonging to $\N_{k-1}, r_{i} \neq r_{j}\
(i,j\in\N_{q})$ such that the $L^{2}(0,\infty)$ spans of the $(k-1-q)$
functions $\{\tilde{\phi}(\cdot);\ s\in\N_{k-1},s \neq r_{i},i\in
\N_{q}\}$ are null.
}
\end{pot}

Now we construct functions $\{\pi_{s}(x;\cdot);s\in \N_{k-p}\}$ defined
on $[0,x], x \in [0,\infty)$ from the set of functions
$\{\phi_{s}(\cdot); s\in \N_{k}, s\neq t_{i},i \in \N_{p}\}$ such that
\setcounter{section}{3}
\begin{equation}
\pi_{s} (x;\cdot) = \sum\limits_{t=1}^{k-p} \alpha_{st}
(x)\phi_{t}(\cdot),\quad (s\in \N_{k-p}, x\in [0,\infty))
\end{equation}
are $(k-p)$ linearly independent functions on $[0,x]$ and
\begin{equation}
\int_{0}^{\pi} \pi_{s} (x;\cdot) \overline{\pi}_{t}(x;\cdot) =
\P_{s}(x)\delta_{st},\quad (s,t\in\N_{k-p}, x\in [0,\infty)),
\end{equation}
where $0 < \P_{1}(x) \leq \P_{2}(x) \leq \ldots \leq \P_{k-p}(x) <
\infty$ and
\begin{equation}
\lim\limits_{x\rightarrow \infty} \P_{s} (x) = \infty.\quad (s\in \N_{k-
p})
\end{equation}
This is achieved by diagonalising the Gram matrix
\begin{equation*}
\Gamma_{x} =
\left[\int_{0}^{x}\phi_{s}(\cdot)\overline{\phi}_{t}(\cdot)\right],\quad
(s,t\in \N_{k-p}, x\in [0,\infty))
\end{equation*}
to diag$(\P_{1}(x), \P_{2}(x),\ldots,\P_{k-p}(x))$ through a unitary
matrix $U_{x} = [\alpha_{st}(x)],\ (s \in \N_{k-p}, t\in\N_{k-p},\ x\in
[0,\infty))$. $\alpha_{st}(x),\ (x\in[0,\infty))$ is the coefficient of
$\phi_{t}(\cdot)$ in the definition (3.17) of $\pi_{s}(x;\cdot)$. Note that
(3.19) is a consequence of the fact that there exists no non-trivial
linear combinations of $\{\phi_{s}; s\in\N_{k-p}\}\in L^{2}(0,\infty)$.
Proceeding along the same lines we construct functions
$\{\sigma_{s}(x;\cdot);\ s\in\N_{k-1-q}\}$ defined on $[0,x], x\in
[0,\infty)$ from the set of functions $\{\tilde{\phi}_{s}(\cdot), \
s\in\N_{k-1}, s\neq r_{i}, i\in \N_{q}\}$ such that
\begin{equation*}
\sigma_{s} (x;\cdot) = \sum\limits_{t=1}^{k-1-q} \beta_{st}
(x)\tilde{\phi}_{t}(\cdot),\quad (s\in\N_{k-1-q}, x\in [0,\infty))
\end{equation*}
and
\begin{equation*}
\int_{0}^{x} \sigma_{s}(x;\cdot) \overline{\sigma}_{t} (x;\cdot) = Q_{s}
(x)\delta_{st},\quad (s,t\in\N_{k-1-q}),
\end{equation*}
where
\begin{equation*}
0<Q_{1}(x) \leq Q_{2}(x) \leq \ldots \leq Q_{k-1-q}(x) < \infty
\end{equation*}
and
\begin{equation*}
\lim\limits_{x\rightarrow\infty} Q_{s}(x) = \infty\ (s \in \N_{k-1-q}).
\end{equation*}
Then we apply Lemma~1.1 to the functions $\{f_{r},g_{r}\}$ given by
\begin{align*}
f_{r}: f_{1},f_{2},\ldots,f_{n+1},\pi_{1}(x;\cdot),\ldots, \pi_{k-p}\
(x;\cdot), \sigma_{1} (x;\cdot),\ldots, \sigma_{k-1-q} (x;\cdot),\\[.2pc]
g_{r}: g_{1},g_{2},\ldots,g_{n+1},\pi_{1}(x;\cdot),\ldots, \pi_{k-p}\
(x;\cdot), \sigma_{1} (x;\cdot),\ldots, \sigma_{k-1-q} (x;\cdot),
\end{align*}
where $f_{r},g_{s}\in \Delta\ \ (r,s\in \N_{n+1})$.

In the resultant determinantal identity evaluated at $x\ (x\in
[0,\infty))$ we divide the $(n + 1 + r)\hbox{th}$ row and column by
$\{P_{r}(x)\}^{1/2}$ and $(n + 1 + k - p + s)\hbox{th}$ row and column by
$\{Q_{s}(x)\}^{1/2}$ for $(r\in \N_{k-p},\ s\in\N_{k-1-q})$. Then we
obtain the following determinant identity given by
\begin{equation*}
\begin{array}{|c|@{\qquad}c@{\qquad}|@{\qquad}c@{\qquad}|}
[f_{r}g_{s}](x)_{r,s\in{\bf N}_{n+1}}
&\displaystyle\frac{[\pi_{r}g_{s}](x)}{P_{r}^{1/2}}
&\displaystyle\frac{[\sigma_{r}g_{s}](x)}{Q_{r}^{1/2}}\\[1pc]\cline{1-3}
 & &\\[-.3pc]
\displaystyle\frac{[f_{r}\pi_{s}](x)}{P_{s}^{1/2}}
&\displaystyle\frac{[\pi_{r}\pi_{s}](x)}{P_{r}^{1/2}P_{s}^{1/2}}
&\displaystyle\frac{[\sigma_{r}\pi_{s}](x)}{Q_{r}^{1/2}P_{s}^{1/2}}\\[1pc]\cline{1-3}
 & &\\[-.3pc]
\displaystyle\frac{[f_{r}\sigma_{s}](x)}{Q_{s}^{1/2}}
&\displaystyle\frac{[\pi_{r}\sigma_{s}](x)}{P_{r}^{1/2}Q_{s}^{1/2}}
&\displaystyle\frac{[\sigma_{r}\sigma_{s}](x)}{Q_{r}^{1/2}Q_{s}^{1/2}}\\[.5pc]
\end{array} = 0.
\end{equation*}

Now we proceed to the limit as $x\rightarrow \infty$. We consider the
limiting values of each of the terms in the above determinant. Note that
\begin{align*}
[\pi_{r}\pi_{r}] (x) &= [\pi_{r}\pi_{r}](0) +
\int_{0}^{x} (\overline{\pi}_{r}M[\pi_{r}] -
\pi_{r}\overline{M}[\pi_{r}])\\
&\quad\ \hbox{(using the Green's formula)}\\
&= [\pi_{r}\pi_{r}](0) +
(\lambda-\overline{\lambda})\int_{0}^{x} \pi_{r}\overline{\pi}_{r}\\
&= [\pi_{r}\pi_{r}](0) + 2i\nu \P_{r}\ (\nu = \hbox{Im}\lambda).
\end{align*}
Therefore,
\begin{equation*}
\frac{[\pi_{r}\pi_{r}](x)}{\P_{r}} = \frac{[\pi_{r}\pi_{r}](0)}{\P_{r}}
+ 2i\nu.
\end{equation*}
Now using (3.19) we get
\begin{equation*}
\lim\limits_{x\rightarrow\infty}\frac{[\pi_{r}\pi_{s}](x)}{\P_{r}} = 2i\nu.
\end{equation*}
Also
\begin{align*}
[\pi_{r}\pi_{s}](x) &= [\pi_{r}\pi_{s}](0) + (\lambda -
\overline{\lambda}) \int_{0}^{x}\pi_{r}\overline{\pi}_{s}\\
&= [\pi_{r}\pi_{s}] (0),\quad \hbox{by (3.18)}.
\end{align*}
Therefore
\begin{equation*}
\lim\limits_{x\rightarrow\infty}\frac{[\pi_{r}\pi_{s}](x)}{P_{r}
^{1/2}P_{s}^{1/2}} = \lim\limits_{x\rightarrow\infty}
\frac{[\pi_{r}\pi_{s}](0)}{P_{r}^{1/2}P_{s}^{1/2}} = 0.
\end{equation*}
Similarly, we get
$\frac{[\sigma_{r}\sigma_{s}](x)}{Q_{r}^{1/2}Q_{s}^{1/2}}$ as
$x\rightarrow \infty$.
\begin{align*}
[f_{r}\pi_{s}](x) &= [f_{r}\pi_{s}](0) + \int_{0}^{x}
(\overline{\pi}_{s}M[f_{r}] - f_{r}\overline{M}[\pi_{s}])\\[.2pc]
&=[f_{r}\pi_{s}](0) + \int_{0}^{x} (\overline{\pi}_{s}M[f_{r}] -
f_{r}\overline{\lambda}\overline{\pi}_{s})\\[.2pc]
&= [f_{r}\pi_{s}](0) + \int_{0}^{x}
\overline{\pi}_{s}(M[f_{r}] - f_{r}\overline{\lambda}),\\[.2pc]
\lim\limits_{x\rightarrow\infty}\frac{[f_{r}\pi_{s}](x)}{P_{s}^{ 1/2}}
&= \lim\limits_{x\rightarrow\infty}\frac{[f_{r}\pi_{s}](0)}{P_{s}^{
1/2}} + \lim\limits_{x\rightarrow\infty} \frac{\int_{0}^{x}
\overline{\pi}_{s}(M[f_{r}] - f_{r}\overline{\lambda})}{P_{s}^{1/2}}\\[.2pc]
&= 0 + \lim\limits_{x\rightarrow\infty}
\frac{\int_{0}^{x}\overline{\pi}_{s}(M[f_{r}]-
f_{r}\overline{\lambda})}{(\int_{0}^{x}|\pi_{s}|^{2})^{1/2}}\\[.2pc]
&= 0.
\end{align*}
(Since $\pi_{s}\in L^{2} (0,x)$, but $\pi_{s}\notin L^{2}(0,\infty)$, we
get the second term also as 0 by Lemma~1.2.) Similarly we can show
that the terms
\begin{equation*}
\frac{[f_{r}\sigma_{s}](x)}{Q_{s}^{1/2}},
\frac{[\pi_{r}g_{s}](x)}{P_{r}^{1/2}},
\frac{[\sigma_{r}g_{s}](x)}{Q_{r}^{1/2}}
\end{equation*}
tend to 0 as $x\rightarrow\infty$. Note that
$[\pi_{r}\sigma_{s}] (x) = [\pi_{r}\sigma_{s}] (0) = 0$.
This follows from the properties of the fundamental solutions, given by
(1.10). Therefore
\begin{equation*}
\lim\limits_{x\rightarrow\infty}
\frac{[\pi_{r}\sigma_{s}](x)}{P_{r}^{1/2}Q_{s}^{1/2}} = 0.
\end{equation*}
Hence taking the limit, the determinant identity becomes
\begin{equation*}
\begin{array}{|c|c@{\quad\ \ }c@{\quad\ \ }c|c@{\quad\ \ }c@{\quad\ \ }c|}
 & & & & & &\\[-.4pc]
[f_{r}g_{s}] \mathop{(\infty)}\limits_{(r,s\in {\bf N}_{n +1})} & &\bigcirc &
& &\bigcirc &\\[1pc]\cline{1-7}
 & & & & & &\\[-.7pc]
 &2i\nu & &\bigcirc & & &\\
\bigcirc & &\ddots & & &\bigcirc &\\
 &\bigcirc & &2i\nu & & &\\[.2pc]\cline{1-7}
 & & & & & &\\[-.7pc]
 & & & &2i\nu & &\bigcirc\\
\bigcirc & &\bigcirc & & &\ddots &\\
 & & & &\bigcirc & &2i\nu\\[.4pc]
\end{array} = 0.
\end{equation*}

Since $\nu \neq 0, \det\limits_{r,s \in {\bf N}_{n+1}} [[f_{r}
g_{s}](\infty)] = 0$ for all $f_{r}, g_{s} \in \Delta$. This completes
the necessity part of the theorem. Now we prove the converse part of the
theorem.

\setcounter{theo}{0}
\begin{rema}
{\rm $\det\limits_{r,s\in {\bf N}_{t}} [[f_{r} g_{s}](\infty)] = 0$ for
any $t < n + 1 \Rightarrow \det\limits_{r,s\in {\bf N}_{n +1}} [[f_{r}
g_{s}](\infty)] = 0$. Hence (2.16) can imply, if at all it is true, that
$p + q \leq n$. Indeed, this is the case which we now prove in the
converse part of the theorem.}
\end{rema}

\begin{suffici}
Assume that (2.16) holds, then we show that $M$ cannot be in the limit
$(k - 1+ p, k + q)$ case at $\infty$ with $p + q > n$. To see this, we
show that $M$ is in the limit $(k -1 + p, k + q)$ case at $\infty$ with
$p + q > n$ contradicting the validity of (2.16).
\end{suffici}

To be specific we assume that $p + q = n + 1\ (p + q > n + 1$ can be
treated along the same lines).

For convenience we define
\begin{align*}
\psi_{r} &:= \psi_{p,r} (x, \lambda),\quad (\lambda \in \C_{+}, r
\in \N_{k -1}),\\[.2pc]
\tilde{\psi}_{r} &:= \psi_{n,r} (x, \overline{\lambda}),\quad (\lambda
\in \C_{+}, r \in \N_{k}),
\end{align*}
where $\psi_{p,r}, \psi_{n,r}$ are defined as in (1.11). Then for a
given $\lambda \in \C_{+}$, in addition to the solutions $\psi_{1},
\psi_{2}, \ldots, \psi_{k-1}$, there exist $p$ solutions which are
linear combinations of $\{\phi_{s} (\cdot, \lambda); s\in \N_{k},
\lambda \in \C_{+}\}$, say
\begin{align}
\xi_{t_{1}} (\cdot) &= \phi_{t_{1}} (\cdot) + \sum\limits_{\substack{s
=1\\ s\neq t_{1}}}^{k} \alpha_{1s} \phi_{s}(\cdot), \quad (t_{1} \in
\N_{k}),\nonumber\\[.1pc]
\xi_{t_{2}} (\cdot) &= \phi_{t_{2}} (\cdot) + \sum\limits_{\substack{s
=1\\ s\neq t_{1},t_{2}}}^{k} \alpha_{2s} \phi_{s}(\cdot), \quad (t_{2} \in
\N_{k}) (t_{2} \neq t_{1}),\nonumber\\
\vdots &\nonumber\\
\xi_{t_{p}} (\cdot) &= \phi_{t_{p}} (\cdot) + \sum\limits_{\substack{s
=1\\ s\neq t_{1},t_{2},\ldots,t_{p}}}^{k} \alpha_{ps} \phi_{s}(\cdot),
\quad (t_{p} \in \N_{k}) (t_{p} \neq t_{1}, t_{2}, \ldots, t_{p-1})
\end{align}
in $L^{2} (0,\infty)$ where $\{\alpha_{is}, s\in \N_{k}, i \in \N_{p}, s
\neq t_{i}\}$ are suitable complex numbers. Note that $\xi_{t_{i}}, i
\in \N_{p}$ belong to $\Delta$.

Similarly for a given $\lambda \in \C_{+}$ in addition to the solutions
$\tilde{\psi}_{1}, \tilde{\psi}_{2}, \ldots, \tilde{\psi}_{k}$ belonging
to $L^{2}(0, \infty)$, there exist $q$ solutions which are linear
combinations of $\{\tilde{\phi}_{s}(\cdot), s \in \N_{k-1}, \lambda \in
\C_{+}\}$ say
\begin{align}
\tilde{\eta}_{r_{1}} (\cdot) &= \tilde{\phi}_{r_{1}} (\cdot) +
\sum\limits_{\substack{s =1\\ s\neq r_{1}}}^{k-1} \beta_{1s}
\tilde{\phi}_{s} (\cdot), \quad (r_{1} \in \N_{k-1}),\nonumber\\[.1pc]
\tilde{\eta}_{r_{2}} (\cdot) &= \tilde{\phi}_{r_{2}} (\cdot) +
\sum\limits_{\substack{s =1\\ s\neq r_{1},r_{2}}}^{k-1} \beta_{1s}
\tilde{\phi}_{s}(\cdot), \quad (r_{2} \in \N_{k-1}) (r_{2} \neq
r_{1}),\nonumber\\
\vdots &\nonumber\\
\tilde{\eta}_{r_{q}} (\cdot) &= \tilde{\phi}_{r_{q}} (\cdot) +
\sum\limits_{\substack{s =1\\ s\neq r_{1},r_{2},\ldots,r_{q}}}^{k-1}
\beta_{1s} \tilde{\phi}_{s}(\cdot), \ \ (r_{q} \in \N_{k-1}) (r_{q}
\neq r_{1}, r_{2}, \ldots, r_{q-1})
\end{align}
in $L^{2} (0,\infty)$ where $\{\beta_{js} : j \in \N_{q}, s \in \N_{k-
1}, s \neq r_{j}\}$ are suitable complex numbers. These solutions
$\tilde{\eta}_{r_{j}} (\cdot), j \in \N_{q}$ belong to $\Delta$.

Now consider the $(n + 1)\times (n + 1)$ determinant, $\det
[[f_{r}g_{s}](\infty)]$ formed by choosing
\begin{align}
f_{r} &: \xi_{t_{1}} (\cdot), \ldots, \xi_{t_{p}} (\cdot),
\psi_{r_{1}} (\cdot), \psi_{r_{2}} (\cdot), \ldots, \psi_{r_{q}}
(\cdot),\nonumber\\[.3pc]
g_{s} &: \tilde{\psi}_{t_{1}} (\cdot), \ldots, \tilde{\psi}_{t_{p}}
(\cdot), \tilde{\eta}_{r_{1}} (\cdot), \tilde{\eta}_{r_{2}} (\cdot),
\ldots, \tilde{\eta}_{r_{q}} (\cdot).
\end{align}
This is given by
\begin{equation}
\begin{array}{|c|c|}
 &\\[-.6pc]
[\xi_{t_{i}} \tilde{\psi}_{t_{j}}]_{i,j\in {\bf N}_{p}} (\infty)
&[\psi_{r_{i}} \tilde{\psi}_{t_{j}}]_{i\in {\bf N}_{q}, j\in {\bf
N}_{p}} (\infty)\\[.5pc]\cline{1-2}
 &\\[-.6pc]
[\xi_{t_{i}} \tilde{\eta}_{r_{j}}]_{i\in {\bf N}_{p}, j\in {\bf N}_{q}}
(\infty) &[\psi_{r_{i}} \tilde{\psi}_{r_{j}}]_{i\in {\bf N}_{p}, j\in
{\bf N}_{q}} (\infty)\\[.5pc]
\end{array}\ .
\end{equation}
We evaluate this determinant (3.23) and show that it is not equal to
zero. From Green's formula (1.6), it follows that all terms in the above
determinant are finite. Further the value of each term of the
determinant is its value at zero. This is because $[\xi\tilde{\eta}](x)$
is independent of $x$. Now we evaluate the determinant as follows. We
consider the two cases separately.

\begin{case}{\rm $t_{i} \neq k$ for every $i$. From the forms of
$\xi_{t_{i}} (\cdot), \tilde{\eta}_{r_{j}} (\cdot)\ (i \in \N_{p}, j \in
\N_{q})$ in (3.20) and (3.21) respectively and the properties of
$\{\theta_{r}\}$ and $\{\phi_{s}\}$ we get
\begin{align*}
\begin{array}{lll}
[\xi_{t_{i}} (x) \tilde{\psi}_{t_{i}} (x)] (\infty) &=\ [\xi_{t_{i}} (x)
\tilde{\psi}_{t_{i}} (x)] (0) &=\ -1\quad (\forall i \in \N_{p}),\\[.55pc]
[\xi_{t_{i}} (x) \tilde{\psi}_{t_{j}} (x)] (\infty) &=\ [\xi_{t_{i}} (x)
\tilde{\psi}_{t_{j}} (x)] (0) &=\ 0\quad (\forall i > j),\\[.55pc]
[\psi_{r_{i}} (x) \tilde{\psi}_{t_{j}} (x)] (\infty) &=\ [\psi_{r_{i}} (x)
\tilde{\psi}_{t_{j}} (x)] (0) &=\ 0\quad (\forall i, j),\\[.55pc]
[\xi_{t_{i}} (x) \tilde{\eta}_{r_{j}} (x)] (\infty) &=\ [\xi_{t_{i}} (x)
\tilde{\eta}_{r_{j}} (x)] (0) &=\ 0\quad (\forall i,j),\\[.55pc]
[\psi_{r_{i}} (x) \tilde{\eta}_{r_{j}} (x)] (\infty) &=\ [\psi_{r_{i}} (x)
\tilde{\eta}_{r_{j}} (x)] (0) &=\ 0\quad (\forall i < j),\\[.55pc]
[\psi_{r_{j}} (x) \tilde{\eta}_{r_{j}} (x)] (\infty) &=\ [\psi_{r_{j}} (x)
\tilde{\eta}_{r_{j}} (x)] (0) &=\ 1\quad (\forall j).
\end{array}
\end{align*}
Therefore in this case we obtain the determinant as
\begin{equation*}
\begin{array}{|cccc|ccccc|}
 & & & & & & & &\\[-.7pc]
-1 & & & & & & & &\\
-\alpha_{1t_{2}} &-1 & &\bigcirc & & & & &\\[.2pc]
-\alpha_{1t_{3}} &-\alpha_{2t_{3}} &-1 & & & & & &\\
\vdots & & & & & &\bigcirc & &\\
-\alpha_{1t_{p}} &-\alpha_{2t_{p}} &\cdots &-1 & & & & &\\[.5pc]\cline{1-9}
 & & & & & & & &\\[-.3pc]
 & & & &1 &\overline{\beta}_{1r_{2}} & &\cdots &\overline{\beta}_{1r_{q}}\\[.3pc]
 &\bigcirc & & & &1 &\overline{\beta}_{2r_{3}} &\cdots &\overline{\beta}_{2r_{q}}\\
 & & & & &\bigcirc & & &\vdots\\
 & & & & & & & &1\\[.3pc]
\end{array}\ .
\end{equation*}}
\end{case}

\begin{case}{\rm $t_{i} = k$ for some $i$. We choose $t_{1} = k$ (there
is no loss of generality in such a choice). Then we get,
\begin{align*}
\begin{array}{lll}
[\xi_{t_{1}} (x) \tilde{\psi}_{t_{1}} (x)] (\infty) &=\ [\xi_{t_{1}} (x)
\tilde{\psi}_{t_{1}} (x)] (0) &=\ i,\ \hbox{the imaginary unit},\\[.55pc]
[\xi_{t_{j}} (x) \tilde{\psi}_{t_{j}} (x)] (\infty) &=\ [\xi_{t_{j}} (x)
\tilde{\psi}_{t_{j}} (x)] (0) &=\ -1\quad (\forall j \geq 2),\\[.55pc]
[\psi_{t_{i}} (x) \tilde{\psi}_{t_{j}} (x)] (\infty) &=\ [\psi_{t_{i}} (x)
\tilde{\psi}_{t_{j}} (x)] (0) &=\ 0\quad (\forall i > j),\\[.55pc]
[\psi_{r_{i}} (x) \tilde{\psi}_{t_{j}} (x)] (\infty) &=\ [\psi_{r_{i}} (x)
\tilde{\psi}_{t_{j}} (x)] (0) &=\ 0\quad (\forall i, j),\\[.55pc]
[\xi_{t_{i}} (x) \tilde{\eta}_{r_{j}} (x)] (\infty) &=\ [\xi_{t_{i}} (x)
\tilde{\eta}_{r_{j}} (x)] (0) &=\ 0\quad (\forall i, j),\\[.55pc]
[\psi_{r_{i}} (x) \tilde{\eta}_{r_{j}} (x)] (\infty) &=\ [\psi_{r_{i}} (x)
\tilde{\eta}_{r_{j}} (x)] (0) &=\ 0\quad (\forall i < j),\\[.55pc]
[\psi_{r_{j}} (x) \tilde{\eta}_{r_{j}} (x)] (\infty) &=\ [\psi_{r_{j}} (x)
\tilde{\eta}_{r_{j}} (x)] (0) &=\ 1\quad (\forall j).
\end{array}
\end{align*}
Therefore in this case the determinant (3.23) takes the form
\begin{equation*}
\begin{array}{|cccc|ccccc|}
 & & & & & & & &\\[-.7pc]
i & & & & & & & &\\
-\alpha_{1t_{2}} &-1 & &\bigcirc & & & & &\\[.2pc]
-\alpha_{1t_{3}} &-\alpha_{2t_{3}} &-1 & & & &\bigcirc & &\\
\vdots & & & & & & & &\\
-\alpha_{1t_{p}} &-\alpha_{2t_{p}} &\cdots &-1 & & & & &\\[.5pc]\cline{1-9}
 & & & & & & & &\\[-.3pc]
 & & & &1 &\overline{\beta}_{1r_{2}} & &\cdots &\overline{\beta}_{1r_{q}}\\[.3pc]
 &\bigcirc & & & &1 &\overline{\beta}_{2r_{3}} &\cdots &\overline{\beta}_{2r_{q}}\\
 & & & & &\bigcirc & & &\vdots\\
 & & & & & & & &1\\[.3pc]
\end{array}\ .
\end{equation*}
Thus in both the cases we get
\begin{equation*}
\det\limits_{r,s\in {\bf N}_{n + 1}} [[f_{r} g_{s}](\infty)] \neq 0,
\end{equation*}
where $\{f_{r}\}$ and $\{g_{s}\}$ are as given in (3.22). This
contradicts hypothesis (2.16). Hence if $\det\limits_{r,s\in {\bf N}_{n + 1}}
[[f_{r} g_{s}](\infty)] = 0$ for all $f_{r}, g_{s} \in \Delta$, then $M$
is in the limit $(k - 1 + p, k + q)$ case at $\infty$ where $p + q \leq
n$.}
\end{case}

\begin{rema}{\rm In the even-order case, the situation of Case~2 does
not arise; this part of the proof follows as in Case~1 (see \cite{1},
ch.~4).}
\end{rema}

Next we prove that $p + q = n$ if (2.16) holds for all $f_{r}, g_{s} \in
\Delta (r, s\in \N_{n + 1})$ but $\det\limits_{r,s\in {\bf N}_{n}} [[f_{r}
g_{s}](\infty)] \neq 0$ for some $f_{r}, g_{s} \in \Delta$. Assume that
$\det\limits_{r,s\in {\bf N}_{n + 1}} [[f_{r} g_{s}](\infty)] = 0$ for all
$f_{r}, g_{s} \in \Delta$. This $\Rightarrow M$ is in the limit $(k - 1
+ p, k + q)$ case at $\infty$, for some $p$ and $q$ such that $p + q
\leq n$.

Further assume that $\det\limits_{r,s\in {\bf N}_{n}} [[f_{r} g_{s}](\infty)]
\neq 0$ for some $f_{r}, g_{s} \in \Delta$. From the necessity part this
$\Rightarrow M$ is in the limit $(k - 1 + p, k + q)$ case at $\infty$
where $p + q \geq n$. Together, it now follows that $M$ is in the limit
$(k - 1 + p, k + q)$ case at $\infty$ where $p + q = n$.

\section{Concluding remarks}

\begin{enumerate}
\renewcommand\labelenumi{\arabic{enumi}.}
\item When $p = 0, q = 0$, we recover the limit-point characterisation
in Theorem~(1.1). Therefore this result is a generalisation of the
characterisation for the limit-point case at $\infty$.

\item Note that, $2k - 1 + n\ (k - 1 + p + k + q)$ gives the gap between
$D (T_{\min})$ and $D(T_{\max})$. More precisely if $(N_{+}, N_{-})$
gives the deficiency indices of $T_{\min}, N_{+} + N_{-}$ gives the
dimension of the quotient space $D(T_{\max})/D(T_{\min})$. Therefore the
above result characterises the gap between $D(T_{\max})$ and
$D(T_{\min})$. However, since $p + q$ and not $p$ and $q$ separately
appear in the theorem the exact break up of the deficiency indices $(k -
1 + p, k + q)$ is not characterised by this result.

\item In the second-order case, this result has been made use of to
obtain sufficient conditions on the coefficients for $M$ to be not in
the limit-point case at $\infty$ \cite{9}.

\item In the third-order case also, this result has been used to obtain
sufficient conditions on the coefficients for $M$ to be not in the
limit-point case at $\infty$ \cite{11}.

\hskip 1pc In the third-order case, the admissible limit classifications of $M$ are
as follows: $M$ is either in the limit (1, 2) case or in the (2, 2)
case, or in the limit (3, 3) case at $\infty$. The above theorem takes
the following forms in these cases:

i. $M$ is in the limit (1, 2) case at $\infty \Leftrightarrow
[fg](\infty) = 0, \quad \forall f, g \in \Delta$.

ii. $M$ is in the limit (2, 2) case at
\begin{equation*}
\hskip -2.5pc \infty \Leftrightarrow \begin{vmatrix}
[f_{1} g_{1}] &[f_{2} g_{2}]\\[.1pc]
[f_{2} g_{1}] &[f_{2} g_{2}]
\end{vmatrix} (\infty) = 0
\end{equation*}
and there exist $f, g \in \Delta$ such that $[fg](\infty) \neq 0$.

iii. $M$ is in the limit (3, 3) case at

\begin{equation}
\hskip -2.5pc \infty \Leftrightarrow \det\limits_{r, s\in {\bf N}_{4}} [[f_{r}
g_{s}](\infty)] = 0,\quad \forall f_{r}, g_{s} \in \Delta
\end{equation}
and there exists $f_{r}, g_{s}\in \Delta$ such that

\begin{equation}
\hskip -2.5pc \det\limits_{r, s\in {\bf N}_{3}} [[f_{r} g_{s}](\infty)] \neq 0.
\end{equation}
Since (4.24) holds irrespective of the limit classifications of $M$
(Lemma~1.1), the required condition reduces to (4.25).

\item We hope that the criteria (i), (ii) and (iii) above may be made
use of to obtain sufficient conditions on the coefficients of $M$ in the
third-order case for $M$ to be in the limit (3,~3) case at $\infty$.

\item Scope for the application of this result to obtain conditions on
the coefficients for specific limit classifications is very much limited
in equations of order higher than 3.

For instance, consider the real, fourth-order case. Here the possible
cases are limit (2,~2), limit (3, 3) and limit (4, 4). Now we can prove
that $M$ is in the limit (4, 4) case at $\infty$, if we can construct
functions $f_{1}, f_{2}, f_{3} \in \Delta$ such that
\begin{equation*}
\hskip -1.25pc \begin{vmatrix}
[f_{1}f_{1}] &[f_{2}f_{1}] &[f_{3}f_{1}]\\[.1pc]
[f_{1}f_{2}] &[f_{2}f_{2}] &[f_{3}f_{2}]\\[.1pc]
[f_{1}f_{3}] &[f_{2}f_{3}] &[f_{3}f_{3}]
\end{vmatrix} (\infty) \neq 0.
\end{equation*}
But this is not an easy task. Also to make use of it in the positive
direction, it is even more difficult. This is so since we have to show
that an $n \times n$ determinant vanishes at the singular point for
every $f_{r}, g_{s} \in \Delta$ for large $n$.

\item If the BVP considered is $My = \lambda wy$ with a positive weight
function $w$ in $[0, \infty)$ we have to set the operators in $L_{w}^{2}
(0, \infty)$, the weighted Hilbert space of integrable square function
with weight $w$ on $[0, \infty)$. The related operators are to be
defined from $\frac{1}{w}M[\cdot]$ acting over functions in $L_{w}^{2} (0,
\infty)$. The sesquilinear form $[fg](x)$ remains the same, but the
functions $f$ and $g$ are now to be chosen in $L_{w}^{2} (0,\infty)$.
The analysis goes through with no significant change.

\item We have discussed the problem on the semi-infinite interval
$[0,\infty)$. However, the analysis goes through without any significant
changes, if it is considered on an interval $[a, b)$ with $-\infty < a <
b \leq \infty$ where $b$ is a singular end-point, i.e., we shall assume
that differential expression is regular on all compact sub-intervals of
$[a, b)$ and is singular at $b$.
\end{enumerate}

\section*{Acknowledgements}

One of the authors (KVA) acknowledges the financial assistance received
from the University of Calicut and the University Grants Commission,
India during the work. AP acknowledges the financial assistance
received from the University of Calicut and the Council of Scientific
and Industrial Research, India during the work.

\end{document}